\documentclass{article}
\usepackage {amssymb}

\newtheorem{theorem}{Theorem}[section]

\newtheorem{lm}[theorem]{Lemma}

\def\irr#1{{\rm  Irr}(#1)}
\def\cd#1{{\rm cd} (#1)}

\title{{\bf Semi-extraspecial Groups}}
\author{ {\bf Mark. L. Lewis}\\[0.3cm]
{\em Department of Mathematical Sciences, Kent State
University,}\\ {\em  Kent, Ohio $44242$, United States of
America}\\
{\em Email}: {\tt lewis@math.kent.edu}  }

\begin{document}

\maketitle

\begin{abstract}
We survey the results regarding semi-extraspecial $p$-groups.  Semi-extraspecial groups can be viewed as generalizations of extraspecial groups.  We present the connections between semi-extraspecial groups and Camina groups and VZ-groups, and give upper bounds on the order of the center and the orders of abelian normal subgroups.  We define ultraspecial groups to be semi-extraspecial groups where the center is as large as possible, and demonstrate a connection between ultraspecial groups that have at least two abelian subgroups whose order is the maximum and semifields.


\end{abstract}

\section{Introduction}

We have two main goals for this paper.  The first is to give an expository account of the known results regarding semi-extraspecial groups.  The second is to connect semifields with a certain class of ultraspecial groups and then show how a number of results from finite geometry apply to these groups.

A nonabelian $p$-group $G$ is {\it special} if $G' = Z(G) = \Phi (G)$.  Furthermore, a group $G$ is {\it extraspecial} if $G$ is a special $p$-group and  $|G'| = |Z (G)| = p$.  Extraspecial groups are central extensions of $Z_p$ by $Z_p^e$, and as their name suggests, they are very special and come in just two types.  We say that a $p$-group $G$ is {\it semi-extraspecial}, if $G$ satisfies the property for every maximal subgroup $N$ of $Z(G)$ that $G/N$ is an extraspecial group.

As far as we can tell, semi-extraspecial groups were first studied by Beisiegel in 1977 \cite{beis}.  The paper \cite{beis} is in German.  One of the other primary places where results regarding semi-extraspecial groups is the paper \cite{ver} by Verardi and concerns a number of the geometric aspects of these groups.  Many of the results regarding semi-extraspecial groups have not appeared in English.  It is important to place these groups into the larger context of Camina groups and VZ-groups. 
 
A group $G$ is a {\it Camina group} if for every element $g \in G \setminus G'$, the conjugacy class of $g$ is $gG'$.  It is in the context of Camina groups that other results about semi-extraspecial groups have been proved.  The main contributions have been by MacDonald, but also by Chillag, Dark, Mann, and Scoppola (see \cite{ChMc}, \cite{DaSc}, \cite{some}, \cite{more}, \cite{mann}, and \cite{MaSc}).  A group $G$ is a VZ-group if every nonlinear irreducible character vanishes off of the center of $G$.  I. e., $\chi (g) = 0$ for every nonlinear character $\chi \in \irr G$ and for every element $g \in G \setminus Z (G)$.  An important source regarding VZ-groups is the paper \cite{extreme} by Fern\'andez-Alcober and Moret\'o, which looks at generalizing several of the key properties of semi-extraspecial groups.

In Section \ref{gens}, we will present a number of equivalent conditions for a group to be a Camina group.  And a similar looking set of equivalent conditions for when a group is a VZ-group.  A group G is a semi-extraspecial group if and only if it is a Camina group and a VZ-group.  Combining the list of conditions, we obtain a list equivalent conditions for a special group to be a semi-extraspecial group.

In Theorem \ref{(Verardi)}, Verardi has proved that there exists an upper bound on the order of abelian subgroups of semi-extraspecial groups.  We see that a particular class of semi-extraspecial groups can be better understood in terms of the number of abelian subgroups whose order equals this upper bound.

The particular class of semi-extraspecial groups that we focus on are the ultraspecial groups.  We will define ultraspecial groups in Section \ref{abel}.  In particular, we are interested in the ultraspecial groups that have at least two abelian subgroups of the largest possible order.  We note that these are not all of the ultraspecial groups.  In particular, Verardi has shown that there exist examples of semi-extraspecial groups with no abelian subgroups of this maximal order and one abelian subgroup of the maximal order.  Verardi also notes that there are many examples known of semi-extraspecial groups with at least two abelian subgroups of maximal order.

We will show that there is a way of constructing a group from the algebraic object called a semifield.  We will will give the definition of semifields in Section \ref{semi}, and we then give the definition of the group associated to each semifield which we will call a {\it semifield group}.  The groups produced by this method from semifields are ultraspecial groups with at least two abelian subgroups.  Semifields were initially studied algebraically as division algebras by Albert and Dickson.  For Albert see: \cite{albert1}, \cite{albert2}, and \cite{albert3} and for Dickson see: \cite{dickson1} and \cite{dickson2}.

The term semifield seems to have originated with Knuth with the papers \cite{knuth1} and \cite{knuth2}.  In fact, \cite{knuth2} contains the results from Knuth's dissertation.  Knuth studied semifields in the context of finite geometries.  Semifields have continued to be an active topic in finite geometries.  For each semifield, one can associate a group .  We will call these groups semifield groups.  In this paper, we review the results regarding semifields that we apply to obtain results regarding semifield groups.  There are many more results known about semifields.  We suggest consulting the expository papers \cite{survey} and \cite{kantor} for a more thorough introduction to semifields.  We believe that there are many other results regarding semifields that can be used to obtain information about semifield groups.

We will see that the construction of these semifield groups has appeared in the finite geometry literature a number of times (see \cite{cronheim}, \cite{hir}, \cite{knst1}, and \cite{roro}).  Our second goal is to present the results from the finite geometry literature regarding these groups.

We would like to thank Professor Alireza Moghaddamfar for his assistance in transcribing the original notes that were the starting point for this paper and for several useful comments while writing this paper.  We would also like to thank Professor James Wilson for helpful discussions while writing the paper and Josh Maglione for computer computations with these groups and many helpful comments while writing this paper.

\section{Overview}

In this paper, all groups are finite $p$-groups for some prime $p$.  This is an expository paper where we present results regarding semi-extraspecial groups and ultraspecial groups, but we will not present proofs.  We will provide references for the proof.

The first time we seriously came to grips with semi-extraspecial groups was when we were writing the paper \cite{VZ} where we determined when two semi-extraspecial groups have the same character table.  We had encountered these groups earlier when we were working on the paper \cite{LMW}, although their involvement in that paper is peripheral as quotients of Camina $p$-groups of nilpotence class $3$.

The genesis of this current paper is several fold.  When presenting the paper \cite{VZ} at a conference, we were asked how many semi-extraspecial groups there are, and then at a later conference where we presenting the paper \cite{Brauer}, some one commented that there were many semi-extraspecial groups.  However, no one seems to be able to quantify how many semi-extraspecial groups exist in ``many.''  In our paper with James Wilson \cite{LeWi}, we give an estimate on the number of quotients of a Heisenberg group, but that paper does not touch the number of ultraspecial groups or the other semi-extraspecial groups that are not quotients of a Heisenberg group.  We will see that this paper gives a way to count one class of ultraspecial groups, but further work is needed to count all ultraspecial groups and in fact to count other semi-extraspecial groups.

Perhaps a more immediate motivation for this paper is to better understand Camina groups of nilpotence class $3$.  It is known that the quotient of a Camina group of nilpotence class $3$ by its center is an ultraspecial group with at least one abelian subgroup of maximal possible order (see Theorem 5.2 of \cite{some} and Theorem 1.3 of \cite{MaSc}).  However, until our recent paper \cite{cents}, all of the examples of Camina $p$-groups that have been constructed have had the Heisenberg group as the quotient modulo the center.  In \cite{cents}, we construct two Camina groups of nilpotence class whose quotients modulo their centers are not isomorphic to the Heisenberg group.   We note that these two groups have quotients modulo the center that are not isomorphic. This raises the question of which ultraspecial groups can occur as the quotient modulo the center of a Camina group of nilpotence class $3$.  Also, in Corollary 6.5 of \cite{cents}, we show that a certain subgroup of Camina groups of nilpotence class $3$ are semi-extraspecial.  With this in mind, we believe that understanding the structure of Camina groups of nilpotence class $3$ will be linked to understanding the structure of semi-extraspecial groups.

We also note that semi-extraspecial groups and ultraspecial groups have appeared recently in a number of classifications or examples of groups with certain properties on their characters and/or their conjugacy classes.  Some of the papers we have in mind are \cite{DMG}, \cite{GGLMNT}, and \cite{KVII}. It seems to us that the understanding of the groups arising in these situations would be enhanced by increasing the knowledge that we have regarding semi-extraspecial groups and ultraspecial groups.

In particular, we hope to encourage new research on these groups.  In this paper, we will focus on abelian subgroups of maximal possible order of ultraspecial groups.  We would believe that results regarding abelian subgroups of maximal possible order for semi-extraspecial groups that are not ultraspecial will deepen our understanding of these groups.  We note that for ultraspecial groups, these abelian subgroups of maximal possible order are centralizers.  We also believe that there would be much benefit to understanding the structure of centralizers in ultraspecial groups that are not abelian.  For semi-extraspecial groups that are not ultraspecial, the centralizers are not abelian.  However there seems to be a close connection between centralizers in these groups and the abelian subgroups.

We also would like to see a better understanding of the quotients of semi-extraspecial groups, and we believe it would be fruitful to determine the automorphism groups of semi-extraspecial and ultraspecial groups.

\section{Extraspecial Groups} \label{extra}

Semi-extraspecial groups can be thought of as a particular type of extension or a generalization of extraspecial groups.  There are many extensions of extraspecial groups that are not semi-extraspecial.  In this section, we review the properties of extraspecial groups.  We highlight the properties that are generalized by semi-extraspecial groups.

Recall that a nonabelian $p$-group $G$ is {\it special} if $G' = Z(G) = \Phi (G)$.  Furthermore, a group $G$ is {\it extraspecial} if $G$ is a special $p$-group and  $|G'| = |Z (G)| = p$.

Here are some of the important facts about extraspecial groups.  These facts are proved in group theory texts such as \cite{hup} or \cite{isagroup}.

\begin{enumerate}
\item[{\rm (1)}] Because $\Phi (G) = G'$, we know that $G/G'$ is elementary abelian; that is, a vector space over ${\rm GF} (p)$.

\item[{\rm (2)}] $|G:G'|$ is a square, so $|G:G'| = p^{2a}$ for some positive integer $a$. (Satz III.13.7 (c) of \cite{hup})

\item[{\rm (3)}] For every prime $p$ and every positive integer $a$, there exists, up to isomorphism, exactly two extraspecial groups of order $p^{2a+1}$.

\item[{\rm (4)}] For every element $x \in G \setminus G'$, the conjugacy class ${\rm cl} (x) = xG'$, so $|C_G (x)| = p^{2a} = |G:G'|$.

\item[{\rm (5)}] If $|G| = p^{2a+1}$, then every noncentral element of $G$ lies in an abelian subgroup of order $p^{a+1}$.

\item[{\rm (6)}] Every irreducible character of $G$ has degree $p^a$ and ``vanishes'' (i.e., is $0$) on $G \setminus Z(G)$. (Theorem 7.5 and Example 7.6 (b) of \cite{hupch})

\item[{\rm (7)}] Every extraspecial group is a central product of extraspecial groups of order $p^3$. (Satz III.13.7 (d) of \cite{hup})
\end{enumerate}

When $p$ is odd, the pair of extraspecial groups of order $p^{2a+1}$ can be distinguished by the fact that one has exponent $p$ and the other has exponent $p^2$.  When $p=2$, it is more complicated to distinguish these groups.  It is well-known that a group of exponent $2$ must be abelian, so we cannot have extraspecial groups of exponent $2$.  It is not difficult to see that that every extraspecial $2$-group will have exponent $4$.  The two extraspecial groups of order $8$ are $Q_8$, the quaternion group of order $8$, and $D_8$, the dihedral group of order $8$.  For any positive integer $a$, the pair of extraspecial groups of order $2^{2a+1}$ can be distinguished by the number of involutions.

Applying (7) above, every extraspecial $2$-group will be a central product of $Q_8$'s and $D_8$'s.  One can show that the central product of two $Q_8$'s is isomorphic to the central product of two $D_8$'s (see the proof of Satz III.13.8 of \cite{hup}).  Hence, a deep theorem shows that the pair of extraspecial groups of order $2^{2a+1}$ can also be distinguished by whether the number of quaternion (or dihedral) factors in the central product is even or odd.

\section{Semi-Extraspecial Groups}

Recall that a $p$-group $G$ is {\it semi-extraspecial}, if $G$ satisfies the property for every maximal subgroup $N$ of $Z(G)$ that $G/N$ is an extraspecial group.  We will often abbreviate semi-extraspecial by s.e.s.  Beisiegel showed that every semi-extraspecial group is a special group (Lemma 1 of \cite{beis}). So if $G$ is s.e.s., then $G' = Z (G) = \Phi (G)$. He also proved that $|G:G'|$ is a square and that $|G'| \leqslant \sqrt{|G:G'|}$ (Satz 1 of \cite{beis}).  In particular, this gives the lemma:

\begin{lm}
Let $G$ be a s.e.s. $p$-group for a prime $p$.  If $|G:G'| = p^{2a}$ and $|G'| = p^b$  for positive integers $a$ and $b$, then $b \leqslant a$.
\end{lm}

In addition, Beisiegel gave a number of examples of semi-extraspecial groups that appear as Sylow subgroups.  In particular, he showed that a Sylow $p$-subgroup of ${\rm SL}_3 (p^a)$ or ${\rm SU}_3 (p^{2a})$ will be s.e.s. of order $p^{3a}$ (see Lemmas 4 and 5 of \cite{beis}).  He also shows that the Suzuki $2$-groups of Types B, C, and D are semi-extraspecial groups (Satz 2 of \cite{beis}).  We say that the group $G$ is the {\it Heisenberg $p$-group of degree $a$} if $G$ is isomorphic to a Sylow $p$-subgroup of ${\rm GL}_3 (p^a)$.  The existence of the Heisenberg groups shows that there exist s.e.s. groups $G$ with $|G:G'| = p^{2a}$ and $|G'| = p^a$ for every prime $p$ and every positive integer $a$.

It is not difficult to prove the following fact.

\begin{lm}
If $G$ is a s.e.s. $p$-group for some prime $p$ and $N < G'$, then $G/N$ is a s.e.s. group.
\end{lm}

When $G$ is the Heisenberg $p$-group of degree $a$ and $b$ is an integer satisfying $1 \leqslant b \leqslant a$, we can find a subgroup $N$ of $G'$ so that $|G':N| = p^b$.  Thus, for every prime $p$ and for all pairs of integers $a, b$ with $1 \leqslant b \leqslant a$, we can find a s.e.s. group $G$ with $|G:G'| = p^a$ and $|G'| = p^b$.

The next result follows from Lemma 2.2 of \cite{some}.  It is also noted at the beginning of Section 2 of \cite{ver}.

\begin{lm}
If $G$ is a s.e.s. $p$-group, then $G/G'$ and $G'$ are elementary abelian $p$-groups.
\end{lm}

\section{Generalizations} \label{gens}

Before we go further, we mention two generalizations of s.e.s. groups.  First, we say that a group $G$ is a {\it Camina group} if for every element $g \in G \setminus G'$, the conjugacy class of $g$ is $gG'$.  Camina groups have been studied in a number of places, namely \cite{ChMc}, \cite{DaSc}, \cite{some}, \cite{more}, and \cite{MaSc}.  It has been shown by Dark and Scoppola that if $G$ is a Camina group, then either $G$ is a Frobenius group whose Frobenius complement is abelian, $G$ is a Frobenius group whose Frobenius complement is isomorphic to the quaternions, or $G$ is a $p$-group.  (An alternate proof of this fact is presented in \cite{fix} and \cite{lewis}.)  Furthermore, MacDonald proved in Theorem 3.1 of \cite{more} that a Camina $2$-group will have nilpotence class $2$, and Dark and Scoppola has proved in \cite{DaSc} when $p$ is odd that a Camina $p$-group will have nilpotence class $2$ or $3$.

There have been a number of equivalent conditions to being a Camina group.  The following has been proved as Proposition 3.1 of \cite{ChMc}.

\begin{theorem}\label{camina}
Let $G$ be a group.  Then the following are equivalent:
\begin{enumerate}
\item $G$ is a Camina group.
\item For every element $g \in G \setminus G'$, $|C_G (g)| = |G:G'|$.
\item For every element $g \in G \setminus G'$ and for every element $z \in G'$, there is an element $y \in G$ so that $[g,y] = z$.
\item Every character $\chi \in \irr G$ vanishes on $G \setminus G'$.
\end{enumerate}
\end{theorem}

The key fact is the following theorem proved by Verardi in Theorem 1.2 of \cite{ver}.

\begin{theorem}
A group $G$ is a s.e.s. $p$-group for some prime $p$ if and only if $G$ is Camina group of nilpotence class $2$.
\end{theorem}

Second, we say that a group $G$ is a VZ-group if every nonlinear irreducible character vanishes off of the center of $G$.  I. e., $\chi (g) = 0$ for every nonlinear character $\chi \in \irr G$ and for every element $g \in G \setminus Z (G)$. If $G$ is a VZ-group, then $G$ is nilpotent of class $2$.  We defined the term VZ-group in \cite{VZ}, although groups satisfying this hypothesis had been earlier studied in a number of places including \cite{extreme} and \cite{KVII}.

To characterize VZ-groups, it is helpful to introduce the concept of isoclinism.  Isoclinism was first introduced by Philip Hall \cite{Hall}.  Two groups $G$ and $H$ are {\it isoclinic} if there exist isomorphisms $\alpha: \frac{G}{Z(G)} \rightarrow \frac{H}{Z(H)}$ and $\beta: G' \rightarrow H'$ such that
$$[\alpha(\overline{g_1}), \alpha(\overline{g_2})] = \beta(\overline{[g_1,g_2]}), \ \ \forall g_1, g_2\in G,$$ where $\overline g = g Z(G)$ for all $g \in G$.  It is not difficult to show that isoclinism is an equivalence relation.  If $G$ and $H$ are isomorphic, then $G$ and $H$ are isoclinic.  On the other hand, if $G$ and $H$ are any two extraspecial groups of the same order, then $G$ and $H$ are isoclinic.  Thus, isoclinism is weaker than isomorphism.

In fact, $G$ and $H$ being isoclinic does not even imply that $|G| = |H|$.  On the other hand, a number of properties such as solvability, nilpotence, and nilpotence class are preserved by isoclinism.  Hall showed that $G$ is always isoclinic to a group $H$ such that $Z(H) \leqslant H'$. For groups of nilpotence class $2$, this implies that every group of nilpotence class $2$ is isoclinic to a special group.

Most of the following result was proved by van der Waall and Kuisch as Theorem 2.4 of \cite{KVII}.  The result as stated is Theorem A of \cite{extreme} by Fern\'andez-Alcober and A.~Moret\'o.

\begin{theorem}\label{VZ}
Let $G$ be a group.  The following conditions are equivalent:
\begin{enumerate}
\item $G$ is a VZ-group.
\item $G$ is isoclinic to a semi-extraspecial $p$-group for some prime $p$.
\item $\cd G = \{ 1, \sqrt {|G:Z(G)|} \}$.
\item For every element $x \in G \setminus Z(G)$, the conjugacy class of $x$ is $xG'$.
\item $Z (G/N) = Z(G)/N$ for every normal subgroup $N$ of $G$ such that $G' \not\leqslant N$.
\end{enumerate}
\end{theorem}

In Theorem B of \cite{extreme}, they prove:

\begin{theorem}\label{normal}
Let $G$ be a group where $|G:Z(G)|$ is a square.  Then $G$ is a VZ-group if and only if every normal subgroup contains $G'$ or is contained in $Z(G)$.
\end{theorem}

Also, it is not difficult to observe that if $G$ is a VZ-group, then $G$ is a semi-extraspecial if and only if $G$ is special.  With this in mind, we combine the equivalent conditions in Theorems \ref{camina}, \ref{normal} and \ref{VZ} specialized to special groups, and we obtain the following equivalent conditions for being semi-extraspecial.


\begin{theorem} \label{equiv}
Let $G$ be a special group. Then, the following are equivalent:
\begin{enumerate}
\item $G$ is semi-extraspecial.
\item For every element $g \in G \setminus G'$, the conjugacy class of $g$ is $gG'$.
\item For every element $g \in G \setminus G'$, we have $|C_G(g)| = |G:G'|$.
\item For every element $g \in G \setminus G'$ and every element $z \in Z (G)$, there exists an element $y \in G$ such that $[g,y] = z$.
\item Every nonlinear character $\chi \in \irr G$ vanishes on $G \setminus G'$.
\item ${\rm cd} (G) = \{ 1, \sqrt {|G:G'|} \}$.
\item For normal subgroup $N$ in $G$, either $G \leqslant N$ or $N \leqslant G'$.
\end{enumerate}
\end{theorem}

Notice that many of the conditions of this generalize a number of the facts listed about extraspecial groups in Section \ref{extra}.

\section{Character Tables}

We came to these groups from character theory.  Recall that the character table of a group is a $k (G) \times k (G)$ complex valued matrix where $k (G)$ is the number of conjugacy classes of $G$ which also equals the number of irreducible characters.  A pair of groups is said to have {\it isomorphic character tables} if the rows and columns of their character tables can be permuted to give identical matrices.  In \cite{VZ}, we determined which VZ-groups have isomorphic character tables.  The condition we found there can be simplified to the following condition for deciding when two s.e.s. groups have isomorphic character tables.

\begin{theorem}
If $G$ and $H$ are s.e.s. groups, then $G$ and $H$ have isomorphic character tables if and only if $|G:G'| = |H:H'|$ and $|G'| = |H'|$.
\end{theorem}

If $G$ is a $p$-group and $C$ is a conjugacy class of $G$, then $C^p = \{ x^p \mid x \in C \}$ is a conjugacy class of $G$.  The map on the conjugacy classes of $G$ given by $C \mapsto C^p$ is called the {\it power map} of $G$. Note that $G$ has exponent $p$ if and only if the power map on $G$ maps every conjugacy class to the conjugacy class of the identity.  We say that $G$ and $H$ form a {\it Brauer pair} if $G$ is not isomorphic to $H$ and $G$ and $H$ have the same character table so that the power maps on the corresponding classes match up.

Let $G$ be a $p$-group.  Set $\Omega_1 (G) = \langle x \in G \mid x^p = 1 \rangle$, and $\mho_1 (G) = \langle x^p \mid x \in G \rangle$.  In \cite{nenciu}, Nenciu determined which pairs of VZ-groups form Brauer pairs.  We then used Nenciu's result to determine which pairs of s.e.s. groups yield Brauer pairs (see \cite{Brauer}).

\begin{theorem}
Let $P$ and $Q$ be non-isomorphic s.e.s. $p$-groups for an odd prime $p$. Then $P$ and $Q$ form a Brauer pair if and only if $|P:P'| = |Q:Q'|$, $|P'| = |Q'|$, and $|\mho_1(P)| = |\mho_1(Q)|$.  In particular, if $P$ and $Q$ have exponent $p$, $|P:P'| = |Q:Q'|$, and $|P'| = |Q'|$, then $P$ and $Q$ form a Brauer pair.
\end{theorem}

Thus, the previous theorem says that s.e.s. groups with exponent $p$ form Brauer pairs if and only if they have the same character tables.  The following result shows that this gives many different Brauer pairs.  Recall that the Heisenberg groups have exponent $p$, so in light of the previous two theorems, these quotients have the same character tables if they have the same size, and so, they form Brauer pairs if they have same size.  This next theorem appeared in \cite{LeWi}.

\begin{theorem}
For every odd prime $p$ and every integer $n \geqslant 12$, there is a Heisenberg group of order $p^{5n/4+O(1)}$ that has $p^{n^2/24+O(n)}$ pairwise nonisomorphic quotients of order $p^n$.
\end{theorem}


In particular, these groups are not distinguishable by their character tables, and thus, it seems that character theory does not have much to say about these groups.  Before we leave character tables, we want to mention the recent result by Casey Wynn in his dissertation.  He studies the super-character theories of semi-extraspecial groups.  He shows that determining the super-character theories of a semi-extraspecial group $G$ can be reduced to determining the super-character theories of the subgroups of $G/G'$ and $G'$, both of which are elementary abelian groups.  See \cite{wynn} which is being published as \cite{LeWy}.

\section{Abelian Subgroups} \label{abel}

We now turn to looking at abelian subgroups of these groups.  We will see that some of these groups can be distinguished by their abelian subgroups.  We begin by stating the following theorem which Verardi proved as Proposition 1.7 and Theorem 1.8 of \cite{ver}.  This theorem gives an upper bound for the order of abelian subgroups of a s.e.s. group.

\begin{theorem}  \label{(Verardi)}
If $G$ is a s.e.s. $p$-group with $|G:G'|=p^{2a}$ and $|G'|=p^b$, then the following are true:
\begin{enumerate}
\item Every abelian subgroup of $G$ has order at most $p^{a+b}$.
\item For every element $g \in G \setminus G'$, we have $p^{b+1} \leqslant |Z(C_G(g))| \leqslant p^{2b}$.
\end{enumerate}
\end{theorem}

Recall from Theorem \ref{equiv} that $|C_G (g)| = p^{2a}$ for all $g \in G \setminus G'$ where we use the notation of Theorem \ref{(Verardi)}.  It follows from part (2) that the only way that $C_G (g)$ can be abelian is if $a = b$.  Thus, we follow Beisiegel in \cite{beis}, and we say that $G$ is {\it ultraspecial} if $G$ is semi-extraspecial and $|G'| = \sqrt {|G:G'|}$.   The next theorem of Verardi is Corollary 5.11 of \cite{ver}.

\begin{theorem}
For each odd prime $p$, the Heisenberg group of degree $p^2$ is the unique ultraspecial group of order $p^6$ and exponent $p$ up to isomorphism.
\end{theorem}

Using the classification of groups of order $p^6$, one can see that there are $p+3$ ultraspecial groups of order $p^6$. The classification of groups of order $p^6$ can be found in \cite{NOV}.  Section 2 of that paper has a nice history of the problem of classifying groups of order $p^6$.  Using the classification, one can also see for each prime $p$ that all of the ultraspecial groups of order $p^6$ are isoclinic (including for $p = 2$).

Next we consider the intersections of the abelian subgroups that have maximal order.  Notice that if $A$ is an abelian subgroup of $G$, then $A Z(G)$ is also an abelian subgroup.  It follows that if $A$ is an abelian subgroup of $G$ of maximal order, then $Z (G) \leqslant A$.  When $G$ is a s.e.s. group, we have that $Z (G) = G'$, and we deduce that every abelian subgroup of $G$ having maximal order contains $G'$ which implies that abelian subgroups of maximal order are all normal in $G$.  This next theorem of Verardi's was proved as Theorem 1.9 of \cite{ver}.

\begin{theorem}
If $G$ is a s.e.s. group with $|G:G'| = p^{2a}$ and $|G'| = p^b$ and $A$ and $B$ are distinct abelian subgroups of $G$ of order $p^{a+b}$ with $|A \cap B| > p^b$, then $p^{2b} \leqslant |A \cap B|\leqslant p^a.$  In particular, if $b > \frac{a}{2}$, then $A \cap B = G'$.
\end{theorem}

In particular, as noted in Corollary 1.10 of \cite {ver}, this theorem implies for ultraspecial groups of order $p^{3a}$ that distinct abelian subgroups of order $p^{2a}$ intersect in the center.

Next, we consider the question of how many abelian subgroups of maximal possible order a s.e.s. can have.  Verardi provides an example of ultraspecial group of order $3^9$ with no abelian subgroups of order $3^6$ on pages 148-149 of \cite{ver}.  Using the computer algebra system Magma \cite{magma}, we have found $4162$ ultraspecial special groups of order $2^9$ with no abelian subgroups of order $2^6$.  One of which is ${\rm SmallGroups} (512, 10477021)$.  Also, using Magma, we have found an ultraspecial groups of orders $5^9$ and $7^9$ (respectively) with no abelian subgroups of order $5^6$ and $7^6$ (repectively).

We believe it is likely that there exist s.e.s groups $G$ with $|G:G'| = p^{2a}$ and $|G'| = p^b$ for every prime $p$ and for all pairs of integers $a$ and $b$ satisfying $a \geqslant b \geqslant 3$.  However, at this time, we do not see any way to prove this fact.  In fact, we do not see any uniform way of producing these groups.  At this point, finding such groups has been a matter of trial and error.

There are many open questions that one can ask about these groups.  One obvious question is what is the smallest order possible for maximal abelian subgroups.  We note that it is easy to find abelian subgroups of order $p^{b+1}$, but it seems unlikely that these will be maximal abelian subgroups.

On pages 149-150 of \cite{ver}, Verardi produces an example of an ultraspecial group of order $7^{9}$ with one abelian subgroup of order $7^{6}$.  Using Magma, we have found many examples of ultraspecial groups of order $3^{12}$ with one abelian subgroup of order $3^8$.  We also found similar examples for many other orders.  In a future paper, we plan to describe how to construct ultraspecial groups of order $p^{3a}$ with one abelian subgroup of order $p^{2a}$ for any prime $p$ and every integer $a \geqslant 3$.  We also will show that this construction can be generalized to other s.e.s. groups.

Finally, in Example 3.9 of \cite{ver}, Verardi gives examples of various s.e.s. groups $G$ with $|G:G'| = p^{2a}$ and $|G'| = p^b$ with at least two abelian subgroups of order $p^{a+b}$ for all primes $p$ and various choices of $a$ and $b$.  Heineken in \cite{hein} gave examples of s.e.s. groups with exactly two abelian subgroups of order $p^{a+b}$.  Verardi has shown that the set of possible values for the number of abelian subgroups of order $p^{a+b}$ is very limited when $b > \frac {a}{2}$.  The following theorem is a combination of Theorems 3.8 and Corollary 5.9 of \cite{ver}.

\begin{theorem}\label{count}
Let $G$ be a s.e.s. group with $|G:G'| = p^{2a}$ and $|G'| = p^b$.  If $b > \frac{a}{2}$ and $G$ has at least $2$ abelian subgroups of order $p^{a+b}$, then the number of abelian subgroups of $G$ of order $p^{a+b}$ has the form $1 + p^h$ where $h$ is an integer satisfying $0 \leqslant h \leqslant a$.  Furthermore, if $G$ is an ultraspecial group and $h > 0$, then $h$ divides $a$.
\end{theorem}

We will use $h (G)$ to denote the quantity $h$ in Theorem \ref{count}.  As noted above, Heineken has produced examples where $h = 0$.  On the other extreme, the Heisenberg groups have $h=a$.  In fact, the following theorem which is Theorem 5.10 in \cite{ver} due to Verardi shows that the number of s.e.s. groups with $h = a$ are very limited.

\begin{theorem}
If $p$ is an odd prime, $G$ is ultraspecial group of order, $h(G) = a$, and $G$ has exponent $p$, then $G$ is isomorphic to the Heisenberg group for $p$ of degree $a$.
\end{theorem}

Now, it seems reasonable to ask for ultraspecial groups of order $p^{3a}$ which values are possible for $h (G)$.  This seems to depend on $p$ and $a$.  We will see when $p = 3$ and $a = 3$ that the possible values for $h (G)$ are $1$ and $3$, and when $p = 3$ and $a = 4$, the possible values for $h (G)$ are $0$, $2$, and $4$.  On the other hand, for $p = 2$ and $a = 3$, the only possible value for $h (G)$ is $3$; for $p = 2$ and $a = 4$, the possible values for $h (G)$ are $0$ and $4$, and for $p = 2$ and $a = 5$, the possible values for $h (G)$ are $0$, $1$, and $5$ (all that are allowed).  Note that the condition that $h (G)$ must divide $a$ is only proved for ultraspecial groups.  It is an open question as to whether other values of $h (G)$ can occur for the other s.e.s. groups covered by that theorem.  We note that Theorem \ref{count} is a consequence of Corollary 3.6 of \cite{ver} which is the following.

\begin{theorem}\label{count2}
Let $G$ be a s.e.s. group with $|G:G'| = p^{2a}$ and $|G'| = p^b$ abelian subgroups $A$ and $B$ so that $G = AB$.  Then the number of abelian subgroups $C$ of $G$ having order $p^{a+b}$ so that $C \cap A = G'$ is $p^h$.
\end{theorem}

In Example 3.9 (c) of \cite{ver}, Verardi gives an example of a s.e.s. group $G$ with $|G:G'| = 3^8$ and $|G'| = 3^2$ with $G = AB$ for abelian subgroups $A$ and $B$ where $h = 6$ which is greater than $a = 4$.  So this shows that the value for $h$ in Theorem \ref{count2} does not obey the upper bound given in Theorem \ref{count}.  A natural open question is determine what values of $h$ are possible in this case.

\section{Semifields} \label{semi}

For the rest of this paper, we will focus on ultraspecial groups of order $p^{3a}$ with at least two abelian subgroups of order $p^{2a}$.  We begin by looking at semifields.  As we stated in the introduction, semifields were initially studied under the name finite division algebras.  This is the name that Verardi uses in \cite{ver}.  Since he does not refer to any of the papers on semifield, we feel it is likely that Verardi was not aware of the literature regarding semifields.  Also, as we mentioned in the introduction, there are two expository articles regarding semifields \cite{survey} and \cite{kantor}. Most of these basic definitions and results can be found there, and we do not offer any more specific reference.

We have already seen ultraspecial groups of order $p^{3a}$ with at least two abelian subgroups of order $p^{2a}$.  Recall that the Heisenberg $p$-group of degree $a$ is  a Sylow $p$-subgroup of ${\rm GL}_3 (p^a)$.  It is known that the upper triangular matrices with $1$'s on the diagonal form a Sylow $p$-subgroup of ${\rm GL}_3 (p^a)$.  Thus, we may represent the Heisenberg group by
$$
\left\{ \left[ \begin{array}{ccc} 1 & a & c \\ 0 & 1 & b \\ 0 & 0 & 1 \end{array} \right] \left| \right. a, b, c \in {\rm GF} (p^a) \right\}
$$
Where ${\rm GF} (p^a)$ is the finite field of order $p^a$.  It is not difficult to see that this group is isomorphic to the group whose set is $\{ (a,b,c) \mid a,b,c \in {\rm GF} (p^a) \}$ with multiplication given by
$$ (a_1, b_1, c_1) \cdot (a_2, b_2, c_2) = ( a_1 + a_2, b_1 + b_2, c_1 + c_2 + a_1  b_2).$$
As we mentioned before, this groups is an ultraspecial group.  It is not difficult to see that $A_1=\{(a,0,c) \mid a, c \in F \}$ and $A_2=\{(0,b,c) \mid b, c \in F \}$ are abelian subgroups of order $p^{2a}$.

Following the literature, we say $(F,+,\ast)$ is a {\it pre-semifield} if $(F,+)$ is an abelian group with at least $2$ elements whose identity is $0$ and $\ast$ is a multiplication that satisfies the distributive laws and $a \ast c = 0$ for $a, c\in F$ implies $a = 0$ or $c = 0$.  (Note we are {\em not} assuming $\ast$ is associative).  We say $F$ is a {\it semifield} if in addition $F$ has an identity $1$.  When $F$ is finite, one can show that $F$ has a vector space structure over some finite field.  This implies that $|F|$ is a power of a prime $p$.  For the remainder of this paper, all (pre-)semifields will be finite.  Note that if $\ast$ is associative, then $F$ will in fact itself be a field.  If $|F|$ is $p$, $p^2$, or $8$, then $*$ must be associative and $(F,+,*)$ will be a field.  We say that $(F,+,*)$ is a {\it proper} semifield  if $*$ is nonassociative.  It has been shown that there exist proper semifields for every prime power $p^n$ such that $n \geqslant 3$ and $p^n \geqslant 16$.  

Let  $(F,+,\ast)$ be a (finite) pre-semifield.  Fix the set $G (F) = \{(a,b,c) \mid a, b, c \in F\}$.  We define the multiplication $\cdot$ on $G(F)$ by
$$ (a_1, b_1, c_1) \cdot (a_2, b_2, c_2) = ( a_1 + a_2, b_1 + b_2, c_1 + c_2 + a_1 \ast b_2).$$
The following theorem has been proved in a number of places.  See Lemma 3 in \cite{beis}, Lemma 2.1 of \cite{hir}, Proposition 2.2 of \cite{knst1}, the Theorem in \cite{roro}, and page 139 of \cite{ver} when $p$ is odd.

\begin{theorem}
If $F$ is a (pre-)semifield, then $G(F)$ is an ultraspecial group of order $|F|^3$ with at least two abelian subgroups of order $|F|^2$.
\end{theorem}

As far as we can tell, semifield groups were first studied by Cronheim in Section 5 of \cite{cronheim}.  Cronheim calls these groups $T$-groups.  Verardi calls these groups $B$-groups in honor of Beisiegel.  Note when $F$ is a finite field, then $G (F)$ is the Heisenberg group for $F$.  Because of this analogy, Knarr and Stroppel call these groups Heisenberg groups, but with wish to reserve the name Heisenberg groups for only the case when $F$ is a field.  Finally, Hiranme calls these groups semifield groups, and that is the name that we shall adopt for these groups.  In particular, we say that $G(F)$ is the {\it semifield group} associated with $F$.

We can identify two of abelian subgroups of $|F|^2$ in $G (F)$.  It is not difficult to see that $A_1=\{(a,0,c) \mid a, c \in F \}$ and $A_2=\{(0,b,c) \mid b, c \in F \}$
are abelian subgroups of order $|F|^2$.  If $p$ is odd, $G(F)$ has exponent $p$ (see Proposition 2.3 (5) of \cite{knst1} or page 139 of \cite{ver}).  For $p = 2$, the elements in $A_1 \cup A_2$ have order $2$ and every element of $G(F)$ outside of $A_1 \cup A_2$ has order $4$ (see Proposition 2.3 (6) of \cite{knst1}).  The following theorem by Verardi (Proposition 3.1 of \cite{ver}) yields the connection between semifields and ultraspecial groups with two abelian subgroups of maximal possible order.  Note that Proposition 3.5 of \cite{hir} gives a similar result.

\begin{theorem}
If $G$ is an ultraspecial $p$-group for an odd prime $p$ with at least two abelian subgroups of order $|G:G'|$ and $G$ has exponent $p$, then there is a semifield $F$ so that $G \cong G(F)$.
\end{theorem}

When $p = 2$, we know that a nonabelian $2$-group cannot have exponent $2$, the hypotheses of the previous theorem could not be met with $p = 2$.  However, if we only require the two specified abelian subgroups to have exponent $2$, then a similar result is obtained. This follows from Proposition 3.5 of \cite{hir}.

\begin{theorem}
If $G$ is an ultraspecial $2$-group with two abelian subgroups, each having order $|G:G'|$ and each having exponent $2$, then there is a semifield $F$ so that $G \cong G(F)$.
\end{theorem}

We note that Hiranime includes a condition that is equivalent to being a semifield group (Lemma 2.2 and Proposition 3.5 of \cite{hir}).

\begin{theorem}
A $p$-group $G$ of order $p^{3a}$ is a semifield group if and only if there exist elementary abelian groups $A$ and $B$ of order $p^{2a}$ so that $a'b' = b'a'$ implies $a' \in A \cap B$ or $b' \in A \cap B$ for all $a' \in A$ and $b' \in B$.
\end{theorem}


Using the Universal Coefficients Theorem (see Chapter 5 of \cite{warfield}), one can show that every s.e.s. $p$-group is isoclinic to a unique s.e.s $p$-group of the same order with exponent $p$ when $p$ is odd.  Thus, when $p$ is odd, one can determine all of the ultraspecial groups $G$ with at least two abelian subgroups of order $|G:G'|$ by determining all the semifield groups and then determining all the possible isoclinisms for each semifield group.

One can also use the Universal Coefficients Theorem to see that every ultraspecial $2$-groups with at least two abelian subgroups of order $|G:G'|$ is isoclinic to an ultraspecial group with two abelian subgroups of order $|G:G'|$ and exponent $2$.  If this is true, again we can find all of the ultraspecial $2$-groups having at least two abelian subgroups of order $|G:G'|$ by computing the semifield groups and then determining all possible isoclinisms.  We note that for ultraspecial $2$-groups that do not have at least two abelian subgroups of order $|G:G'|$ that it is not so clear how one should determine representatives from each isoclinism class.

\section{More on Semifield Groups}

We have seen that every ultraspecial group $G$ with at least two abelian subgroups is isoclinic to a unique (up to isomorphism) semifield group for some semifield.  We now want to determine when two semifield groups are isomorphic.

We say two (pre-)semifields $(F_1,+,\ast_1)$ and $(F_2,+,\ast_2)$ are {\it isotopic} if there exist additive isomorphisms $\alpha, \beta, \gamma: F_1 \rightarrow F_2$ that satisfy
$$\gamma (a \ast_1 b) = \alpha (a) \ast_2 \beta (b), \ \ \ \ \forall a, b \in F_1.$$
One can show that isotopism is an equivalence relation on (pre-)semifields.  It is well-known in the semifield literature that every pre-semifield is isotopic to some semifield.  We note that it is usually easier to describe a pre-semifield for each isotopism class of semifields; so often it is pre-semifields that are described.  For example, Albert's ``twisted semifields'' are defined as follows.  Let $F$ be a finite field.  For an element $j \in F$ and nontrivial automorphisms $\alpha, \beta$ define $\ast$ by $x \ast y = xy + j \alpha (x) \beta (y)$ (see (3.3) of \cite{kantor}).  Then $(F,+,\ast)$ will be a pre-semifield.

We now begin to see how isotopism of the semifields yields isomorphic semifield groups.  The following result has been proved by a number of authors (see Lemma 2.4 of \cite{hir}, Proposition 3.2 (1) of \cite{knst1}, Proposition A of \cite{roro}, and Lemma 3.13 of \cite{ver}).

\begin{theorem}
If $F_1$ and $F_2$ are isotopic (pre-)semifields, then $G (F_1)$ and $G (F_2)$ are isomorphic groups.
\end{theorem}


If $|F|$ is $p$, $p^2$ or $8$, then $F$ must be isotopic to the field of that order.  (This has been proved by Menichetti.)  Thus, for these orders, there is only one isotopism class of (pre)-semifield.  Equivalently, this says that there is a unique semifield group of order $8^3$.  On the other hand, it is known for every prime power $p^a$ with either $p$ odd and $a \geqslant 3$ or $p = 2$ and $a \geqslant 4$ that there exists some semifield that is {\em not} isotopic to the field of that order.  In Proposition B of \cite{roro}, Rocco and Rocha prove the following:

\begin{theorem} \label{proper}
Let $G$ be the Heisenberg group of degree $a$ for the prime $p$.  If $F$ is a proper semifield of order $p^a$, then $G(F)$ is not isomorphic to $G$.
\end{theorem}

Note that Theorem \ref{proper} implies that for every prime $p$ and integer $a$ satisfying either $p$ is odd and $a \geqslant 3$ or $p = 2$ and $a \geqslant 4$, then there exists a semifield group of order $p^{3a}$ that is not isomorphic to the Heisenberg group.

If $F$ is a (pre)-semifield, then we define $F^{\rm op}$ by $a \ast^{\rm op} b = b \ast a$. Obviously, if $F$ is commutative, then $F^{\rm op} = F$.  It is possible to have $F$ isotopic to $F^{\rm op}$ when $F$ is not isotopic to a commutative semifield.  Both Hiranime (in Theorem 5.1 of \cite{hir}) and Knarr and Stroppel (Proposition 3.2 (b) of \cite{knst1}) have proved:

\begin{theorem}
If $F$ is a semifield, then $G (F) \cong G (F^{\rm op})$.
\end{theorem}

We say that $(F_1,+,\ast_1)$ and $(F_2,+,\ast_2)$ are {\it anti-isotopic} if there exist (additive) isomorphisms $\alpha, \beta, \gamma: F_1 \rightarrow F_2$  such that
$$\gamma (a \ast_1 b) = \beta (b) \ast_2 \alpha(a), \ \ \ \ \forall a,b \in F_1.$$
It is not difficult to see that $F_1$ and $F_2$ are anti-isotopic if and only if $F_1$ and $F_2^{\rm op}$
are isotopic.  We have seen that if $F_1$ and $F_2$ are either isotopic or anti-isotopic, then $G (F_1)$ and $G (F_2)$ are isomorphic.  In fact, these are the only times that these groups can be isomorphic.  This is the content of the next theorem.  This is proved by Hiranime as Theorem 5.1 of \cite{hir} and Knarr and Stroppel as Theorem 6.6 of \cite{knst1}.

\begin{theorem} \label{isom} 
If $F_1$ and $F_2$ are (pre-)semifields, then $G(F_1)$ and $G(F_2)$ are isomorphic if and only if $F_1$ and $F_2$ are isotopic or anti-isotopic.
\end{theorem}

We can define an equivalence relation on semifields by saying that two semifields are equivalent if they are isotopic or anti-isotopic.  Theorem \ref{isom} shows that the isomorphism classes of the semifield groups are in bijection to the equivalence classes of this equivalence relation.  To emphasize this point, the number of semifield groups of order $p^{3a}$ equals the number of unique semifields up to isotopism or anti-isotopism.

As we stated above, the field of order $8$ is the only semifield of order $8$, so the Heisenberg group is the unique semifield group of order $2^9 = 512$.  Kleinfeld in \cite{klein} and Knuth in \cite{knuth1} computed the semifields of order $16$.  There are two isotopism classes of proper semifields of order $16$ and the two classes are not anti-isotopic.  Thus, there exist three semifield groups of order $2^{12}$ up to isomorphism.  Knuth \cite{knuth1} and Walker \cite{walker} determined the semifields of order $32$.  There are six isotopism classes of semfields of order $32$, and of these, two pairs are anti-isotopic.  Thus, there exist four semifield groups of order $2^{15}$.  Finally, using a computer search, R\'ua, Combarro, and Ranilla \cite{RCR1} found all the semifields of order $64$.  They have determined that there are $332$ isotopism classes of semifields of order $64$.  The number of equivalence classes under isotopsim and anti-isotopism is $184$, so the number of semifield groups of order $2^{18}$ is $184$.

The semifields of order $27$ were initially computed by Dickson \cite{dickson1}.  It is known that all of the proper semifields of order $27$ are isotopic, so there are two semifield groups of order $3^9$.  The semifields of order $81$ were computed by Dempwolff \cite{demp} using the computer.  There are $27$ isotopism classes of semifields of order $81$ and once anti-isotopisms are included, there are $19$ equivalence classes under isotopism and anti-isotopism.  It follows that there are $19$ semifield groups of order $3^{12}$.   R\'ua, Combarro, and Ranilla \cite{RCR2} have also used the computer to find the semifields of order $243$.  There are $23$ isotopism classes of semifields of order $243$ and accounting for both isotopism and anti-isotopism, there are $15$ equivalence classes; so the number of semifield groups of order $3^{15}$ is $15$.

The semifields of order $5^3$ were also computed by Dickson \cite{dickson1} with mistakes corrected by Long \cite{long}.  There are four isotopism classes of semifields of order $5^3$.  It is possible that there is one pair of anti-isotopic semifields that are not isotopic.  Thus, we have not determined whether there are three or four semifield groups of order $5^9$.  In \cite{CRR2}, they give enough information to compute the semifields of order $5^4$; however, they have not computed the isotopism classes of the semifields of order $5^4$.  The general case of semifields of order $p^3$ was initially studied by Albert in \cite{albert1}.  Menchetti proved a conjecture of Kaplansky regarding the full classification of these semifields in \cite{mench}.  Using this classification, we should be able to determine the semifield groups of order $p^9$.  Finally, we mention that the semifields of order $7^4$ have also been classified \cite{CRR}, so we should also be able to determine the semifield groups of order $7^{12}$.  We summarize these results in Table \ref{table 1}.

{\begin{table}[ht]  
\begin{tabular}{|l|c|c|c|c|}
\hline
$|F|$ & $\#$ isotopism classes & $|G(F)|$ & $\#$ $G(F)$ & $\#$ commutative semifields \\
\hline
$2^3$ & 1 & $2^9$ & 1 & 1 \\
$2^4$ & 3 & $2^{12}$ & 3 & 1 \\
$2^5$ & 6 & $2^{15}$ & 4 & 2 \\
$2^6$ & 332 & $2^{18}$ & 184 & 2 \\
$2^7$ & ? & $2^{21}$ & ? & 2 \\
\hline
$3^3$ & 2 & $3^9$ & 2 & 2 \\
$3^4$ & 27 & $3^{12}$ & 19 & 2 \\
$3^5$ & 23 & $3^{15}$ & 15 & 7 \\
\hline
$5^3$ & 4 & $5^9$ & ? & 2 \\
$5^4$ & ? & $5^{12}$ & ? & ? \\
\hline
$7^3$ & ? & $7^9$ & ? & ? \\
$7^4$ & 356 & $7^{12}$ & 227 & 2 \\
\hline
\end{tabular}
\caption{Number of semifield groups}\label{table 1}
\end{table}}

In this table, we consider possible values of $|F|$.  For each value of $|F|$, we list the number of isotopism classes of semifields of order $|F|$.  We list $|G (F)|$ (which is $|F|^3$).  We list the number of semifield groups of order $|G (F)|$ which equals the number of equivalence classes under isotopism and anti-isotopism and the number of isotopism classes that contain a commutative semifields which equals the number of semifield groups having more than two abelian subgroups of order $|F|^2$.

\section{Commutative Semifields and Seminuclei}

We now see that if $G (F)$ is a semifield group, then the number of abelian subgroups of order $|F|^2$ can be determined by the structure of $F$.  This next result is proved by Verardi  (Theorem 3.14 of \cite{ver}), Hiranime (Proposition 4.2 (i) of \cite{hir}), and Knarr and Stroppel (Lemma 4.3 of \cite{knst1}).

\begin{theorem}
Let $F$ be a semifield.  Then $G (F)$ has more than two abelian subgroups of order $|F|^2$ if and only if $F$ is isotopic to a commutative semifield.
\end{theorem}

In particular, $F$ is not isotopic to a commutative semifield if and only if $G(F)$ has exactly two abelian subgroups of order $|F|^2$.  We summarize the number of isotopism classes of commutative semifields in Table \ref{table 1}.  For the most part, the references for the results in that table are the same as the references for the last couple paragraphs of the previous section, although we also need to refer to \cite{RC}.

When $F$ is a commutative semifield, we can determine the number of abelian subgroups of $G(F)$ of order $|F|^2$.  We define
$${\rm Mid}(F) = \{ z \in F \mid  x \ast (z \ast y) =  (x \ast z) \ast y, \forall x, y \in F \},$$
which is called the {\it middle seminucleus} of $F$.  When $F$ is finite, ${\rm Mid}(F)$ is a field. The following theorem is proved by Verardi (Corollary 5.9 of \cite{ver}), Hiranime (Proposition 4.3 (ii) of \cite{hir}), and Knarr and Stroppel (Lemma 4.3 of \cite{knst1}).

\begin{theorem}
If $F$ is a commutative semifield, then $|{\rm Mid}(F)| = p^h$ where $1+p^h$ is the number of abelian subgroups of order $|F|^2$ in $G(F)$.
\end{theorem}

It is not difficult to see that $F$ is a vector space over ${\rm Mid} (F)$.  This implies when $|F|=p^a$, that $h$ divides $a$ where $|{\rm Mid}(F)| = p^h$.  Recall that we knew from before that $h \leqslant a$.  We note that the Heisenberg groups are the only semfield groups of order $2^9$ and $2^{12}$ that have more than two abelian subgroups of the maximal possible order.  For $2^{15}$, $2^{18}$, $2^{21}$, $3^9$, $3^{12}$, $5^9$, and $7^{12}$ there is for each order, one semifield group other than the Heisenberg group that has more than two abelian subgroups maximal possible order.  For $2^{15}$ and $2^{21}$, we see from Theorem \ref{count} that the semifield group other than the Heisenberg group with more than two abelian subgroups of maximal order will have exactly three abelian subgroups of maximal order.  The non-Heisenberg group of order $2^{18}$ having more that two abelian subgroups of maximal order will have five abelian subgroups of maximal order.

In a similar vein, the non-Heisenberg group of order $3^9$ and $3^{15}$ having more than two abelian subgroups of maximal order have four abelian subgroups of maximal order.  The non-Heisenberg group of order $3^{12}$ with more than two abelian subgroups of maximal order has ten abelian subgroups of maximal order.  Finally, the non-Heisenberg group of $5^9$ with more than two abelian subgroups of maximal order has six abelian subgroups of maximal order.

Define
$${\rm R}(F) = \{ z \in F \mid (x \ast y) \ast z = x \ast (y \ast z), \forall x, y \in F\}.$$  This is the {\it right seminucleus} of $F$.  When $F$ is finite, $R(F)$ is also a field and $F$ is a vector space over ${\rm R} (F)$.  The following result is due to Verardi (Theorem 5.12 of \cite{ver}).

\begin{theorem}
Suppose $F$ is a commutative semifield of order $p^a$, $p$ is an odd prime, and $|{\rm R} (F)| = p^r$.  For every element $v \in G(F) \setminus G(F)'$, set $Z(v) = Z(C_G(v))$.
\begin{itemize}
\item[{\rm (a)}] If $|Z(v)| = p^{a+k}$, then $r$ divides $k$.
\item[{\rm (b)}] If $v_1, \dots, v_m$ in $G (F) \setminus G (F)'$ are chosen so that
$$G (F) \setminus G (F)' = \bigcup_{i=1}^{m} (Z(v_i) \setminus G(F)'),$$ then $m\equiv 1 \pmod {p^r}$.
\end{itemize}
\end{theorem}

Note that the left seminucleus can be defined analogously, and that a similar theorem will hold with the left seminucleus replacing the right seminucleus.  This shows that some information regarding the group structure of $G(F)$ can be obtained from the semifield structure of $F$ when $F$ is commutative.  It seems likely to us that there are other results regarding the structure of semifields that can be applied to obtain group structure of $G(F)$.  We note that there has been much more work studying commutative semifield than general semifields in the literature, and we expect that many of the results regarding commutative semifields can be applied to semifield groups.  Finally, we would like to see results relating the semifield structure of $F$ with the group structure of $G (F)$ when $F$ is not commutative.





\begin{thebibliography}{99}

\bibitem{albert1} A.~A.~Albert, On nonassociative division algebras, {\it Trans. Amer. Math. Soc.} {\bf 72}, (1952), 296-309.

\bibitem{albert2} A.~A.~Albert, Finite noncommutative division algebras, {\it Proc. Amer. Math. Soc.} {\bf 9} (1958), 928-932.

\bibitem{albert3} A.~A.~Albert, Isotopy for generalized twisted fields, {\it An. Acad. Brasil. Ci.} {\bf 33} (1961), 265-275.

\bibitem{beis} B.~Beisiegel, Semi-extraspezielle $p$-Gruppen, {\it Math. Z.} {\bf 156} (1977), 247-254.

\bibitem{magma} W.~Bosma, J.~Cannon, and C.~Playoust. The Magma algebra system I: The user language, {\it J. Symbolic Comput.} {\bf 24} (1997), 235-265.

\bibitem{ChMc} D.~Chillag and I.~D.~MacDonald, generalized Frobenius groups, {\it Israel J. Math.} {\bf 47} (1984), 111-122.

\bibitem{CRR2} E.~F.~Combarro , I.~F.~R\'ua, and J.~Ranilla, New advances in the computational exploration of semifields, {\it Int. J. Comput. Math.} {\bf 88} (2011), 1990-2000.

\bibitem{CRR} E.~F.~Combarro , I.~F.~R\'ua, and J.~Ranilla, Finite semifields with $7^4$ elements, {\it Int. J. Comput. Math.}, {\bf 89} (2012), 1865-1878.

\bibitem{survey} M.~Cordero and G.~P.~Wene, A survey of finite semifields, {\it Discrete Math.} {\bf 208/209} (1999), 125-137.

\bibitem{cronheim} A.~Cronheim, $T$-groups and their geometry, {\it Illinois J. Math.} {\bf 9} (1965), 1-30.

\bibitem{DaSc} R.~Dark and C.~M.~Scoppola, On Camina groups of prime power order, {\it J. Algebra} {\bf 181} (1996), 787-802.

\bibitem{demp} U.~Dempwolff, Semifield planes of order $81$, {\it J. geom.} {\bf 89} (2008), 1-16.

\bibitem{dickson1} L.~E.~Dickson, Linear algebras in which division is always uniquely possible, {\it Trans. Amer. Math. Soc.} {\bf 7} (1906), 370-390.

\bibitem{dickson2} L.~E.~Dickson, On commutative linear algebras in which division is always uniquely possible, {\it Trans. Amer. Math. Soc.} {\bf 7} (1906), 514-522.




\bibitem{DMG} S.~Dolfi, A.~Moret\'o, and G.~Navarro, The groups with exactly one class of size a multiple of $p$, {\it J. Group Theory} {\bf 12} (2009), 219-234.

\bibitem{extreme} G.~A.~Fern\'andez-Alcober and A.~Moret\'o. Groups with two extreme character degrees and their normal subgroups, {\it Trans. Amer. Math. Soc.} {\bf 353} (2001), 2171-2192.

\bibitem{GGLMNT} D.~Goldstein, R.~M.~Guralnick, M.~L.~Lewis, A.~Moret\'o, G.~Navarro,  and P.~H.~Tiep, Groups with exactly one irreducible character of degree divisible by $p$, {\it Algebra Number Theory} {\bf 8} (2014), 397-428.

\bibitem{Hall} P.~Hall, The classification of prime-power groups, {\it J. Reine Angew. Math.} {\bf 182} (1940), 130-141.

\bibitem{hein} H. Heineken,  Nilpotente Gruppen, deren s\"amtliche Normalteiler charakteristisch sind, {\it Arch. Math. (Basel)} {\bf 33} (1979/80), 497-503.

\bibitem{hir} Y.~Hiramine, Automorphisms of $p$-groups of semifield type, {\it Osaka J. Math.} {\bf 20} (1983), 735-746.

\bibitem{hup} B.~Huppert, ``Endliche Gruppen I,'' Springer-Verlag, Berlin-New York 1967.

\bibitem{hupch} B.~Huppert, ``Character Theory of Finite Groups,'' Walter de Gruyter \& Co., Berlin, 1998.


\bibitem{isagroup} I.~M.~Isaacs, ``Finite Group Theory,'' American Mathematical Society, Providence, RI, 2008.


\bibitem{fix} I.~M.~Isaacs, and M.~L.~Lewis, Camina $p$-groups that are generalized Frobenius complements, {\it Arch. Math. (Basel)} {\bf 104} (2015), 401-405.

\bibitem{kantor} W.~M.~Kantor, Finite semifields in {\it Finite geometries, groups, and computation}, 103-114, Walter de Gruyter GmbH $\&$ Co. KG, Berlin, 2006.

\bibitem{klein} E.~Kleinfeld, Techniques for enumerating Veblen-Wedderburn systems, {\it J. Assoc. Comput. Mach.} {\bf 7} (1960) 330-337.

\bibitem{knst1} N.~Knarr and M.~J.~Stroppel, Heisenberg groups, semifields, and translation planes, {\it Beitr. Algebra Geom.} {\bf 56} (2015), 115-127.


\bibitem{knuth1} D.~E.~Knuth, Finite semifields and projective planes, {\it J. Algebra} {\bf 2} (1965), 182-217.

\bibitem{knuth2} D.~E.~Knuth, A class of projective planes, {\it Trans. Amer. Math. Soc.} {\bf 115} (1965), 541-549.




\bibitem{VZ}  M.~L.~Lewis, Character tables of groups where all nonlinear irreducible characters vanish off the center in {\it Ischia group theory 2008}, 174-182, World Sci. Publ., Hackensack, NJ, 2009.

\bibitem{Brauer} M.~L.~Lewis, Brauer pairs of Camina $p$-groups of nilpotence class $2$, {\it Arch. Math. (Basel)} {\bf 92} (2009), 95-98.

\bibitem{lewis} M.~L.~Lewis, Classifying Camina groups: a theorem of Dark and Scoppola, {\it Rocky Mountain J. Math.} {\bf 44} (2014),  591-597.  Erratum on "Classifying Camina groups: a theorem of Dark and Scoppola'' [MR3240515], {\it Rocky Mountain J. Math.} {\bf 45} (2015), 273.

\bibitem{cents} M.~L.~Lewis, Centralizers of Camina groups with nilpotence class $3$,  Submitted.

\bibitem{LMW} M.~L.~Lewis, A.~Moret\'o, and T.~R.~Wolf, Non-divisibility among character degrees, {\it J. Group Theory} {\bf 8} (2005), 561-588.

\bibitem{LeWi} M.~L.~Lewis and J.~B.~Wilson, Isomorphism in expanding families of indistinguishable groups, {\it Groups Complex. Cryptol.} {\bf 4} (2012), 73-110.

\bibitem{LeWy} M.~L.~Lewis and C.~Wynn, Supercharacter theories of semiextraspecial $p$-groups and Frobenius groups, Submitted.

\bibitem{long} F.~W.~Long, Corrections to Dickson's table of three dimensional division algebras over $F_5$, {\it Math. Comp.} {\bf 31} (1977), 1031-1033.

\bibitem{some} I.~D.~Macdonald, Some $p$-groups of Frobenius and extraspecial type, {\it Israel J. Math.} {\bf 40} (1981), 350-364.

\bibitem{more} I.~D.~Macdonald, More on $p$-groups of Frobenius type, {\it Israel J. Math.} {\bf 56} (1986), 335-344.

\bibitem{mann} A.~Mann, Some finite groups with large conjugacy classes, {\it Israel J. Math.} {\bf 71} (1990), 55-63.

\bibitem{MaSc} A.~Mann, and C.~M.~Scoppola, On $p$-groups of Frobenius type, {\it Arch. der Math.} {\bf 56} (1991), 320-332.

\bibitem{mench} G.~Menichetti, On a Kaplansky conjecture concerning three-dimensional division algebras over a finite field, {\it J. Algebra} {\bf 47} (1977), 400-410.


\bibitem{nenciu} A.~Nenciu, Brauer pairs of VZ-groups, {\it J. Algebra Appl.} {\bf 7} (2008), 663-670.

\bibitem{NOV} M.~F.~Newman, E.~A~O'Brien, and M.~R.~Vaughan-Lee, Groups and nilpotent Lie rings whose order is the sixth power of a prime, {\it J. Algebra} {\bf 278} (2004), 383-401.

\bibitem{roro} N.~R.~Rocco and J.~S.~Rocha, A note on finite semifields and certain $p$-groups of class $2$, {\it Discrete Math.} {\bf 275} (2004), 355-362.

\bibitem{RCR1} I.~F.~R\'ua, E.~F.~Combarro, and J.~Ranilla, Classification of semifields of order $64$, {\it J. Algebra} {\bf 322} (2009), 4011-4029.

\bibitem{RC} I.~F.~R\'ua and E.~F.~Combarro, Commutative semifields of order $3^5$, {\it Comm. Algebra} {\bf 40} (2012), 988-996.

\bibitem{RCR2} I.~F.~R\'ua, E.~F.~Combarro, and J.~Ranilla, Determination of division algebras with $243$ elements, {\it Finite Fields Appl.} {\bf 18} (2012), 1148-1155.

\bibitem{saeidi} A.~Saeidi, Finite $2$-groups in which distinct nonlinear irreducible characters have distinct kernels, {\it Quaest. Math.} {\bf 39} (2016), 523-530.

\bibitem{KVII} R.~W.~van der Waall and E.~B.~Kuisch, Homogeneous character induction II, {\it J. Algebra} {\bf 170} (1994), 584-595.

\bibitem{ver} L.~Verardi, Gruppi semiextraseciali di esponente $p$, {\it Ann. Mat. Pura Appl.} {\bf 148} (1987), 131-171.

\bibitem{walker} R.~J.~Walker, Determination of division algebras with $32$ elements, {\it Proc. Sympos. Appl. Math.}, {\bf 15}, (1963), 83-85.

\bibitem{warfield} R.~B.~Warfield, Jr., ``Nilpotent groups,'' Lecture Notes in Mathematics, Vol. 513, Springer-Verlag, Berlin-New York, 1976.

\bibitem{wynn} C.~W.~Wynn, ``Supercharacter theories of Camina pairs,''  Ph.D. dissertation, Kent State University, 2017.
\end{thebibliography}
\end{document}